\documentclass[10pt,twoside]{amsart}

\usepackage{hyperref,amsmath,amsfonts,mathrsfs}
\usepackage{pgf}
\usepackage{tikz-cd}
\calclayout

\usepackage{scalerel}

\newcommand\reallywidehat[1]{\arraycolsep=0pt\relax%
\begin{array}{c}
\stretchto{
  \scaleto{
    \scalerel*[\widthof{\ensuremath{#1}}]{\kern-.5pt\bigwedge\kern-.5pt}
    {\rule[-\textheight/2]{1ex}{\textheight}} 
  }{\textheight} %
}{0.5ex}\\           
#1\\                 
\rule{-1ex}{0ex}
\end{array}
}

\def\pn{\par\noindent}

\def\cen{\centerline}

\begin{document}

\title
{On combinatorial aspects arising from modules over commutative rings}

\author { Rameez Raja }
\maketitle

\pn {\cen {\small \it School of Mathematics, Harish-Chandra Research Institute, HBNI,} 
\pn {\cen {\small \it Chhatnag Road, Jhunsi, Allahabad 211019, India}}
\pn {\cen {\it rameezraja@hri.res.in}}

\begin{abstract}
Let $R$ be a commutative ring with unity, $M$ be an unitary $R$-module and $\Gamma$ be a simple graph. This research article is an interplay of combinatorial and algebraic properties of $M$. We show a combinatorial object completely determines an algebraic object and characterize all finite abelian groups. We discuss the correspondence between essential ideals of $R$, submodules of $M$ and vertices of graphs arising from $M$. We examine various types of equivalence relations on objects $\widehat{A_{f}(M)}$, $\widehat{A_{s}(M)}$ and $\widehat{A_{t}(M)}$, where $\widehat{A_{f}(M)}$ is an object of full-annihilators, $\widehat{A_{s}(M)}$ is an object of semi-annihilators and $\widehat{A_{t}(M)}$ is an object of star-annihilators in  $M$. We study essential ideals corresponding to elements of an object $\widehat{A_{f}(M)}$ over hereditary and regular rings. Further, we study isomorphism of annihilating graphs arising from $M$ and tensor product $M \otimes_{R}T^{-1}R$, where $T = R \backslash C(M)$, where $C(M) = \{ r\in R$ : $rm = 0~ for~ some~ 0\neq m \in M$\}, and show that $ann_f(\Gamma(M \otimes_R T^{-1}R)) \cong ann_f (\Gamma(M))$ for every $R$-module $M$.

\vskip 0.4 true cm



\end{abstract}


\medskip

\section{\bf Introduction}

The subject of algebraic combinatorics deals with the study of combinatorial structures arising in an algebraic context, or applying algebraic techniques to a combinatorial problem. One of the areas in algebraic combinatorics introduced by Beck \cite{Bk} is to study the interplay between graph theoretical and algebraic properties of an algebraic structure. Continuing the concept of associating a graph to an algebraic structure another combinatorial approach of studying commutative rings was given by Anderson and Livingston in \cite{AdLs}. They associated a simple graph to a commutative ring $R$ with unity called a zero-divisor graph denoted by $\Gamma(R)$ with vertices as $Z^{*}(R) = Z(R)\backslash \{0\}$, where $Z(R)$ is the set of zero-divisors of $R$. Two distinct vertices $x, y\in Z^{*}(R)$ of $\Gamma(R)$ are adjacent if and only if $xy = 0$. The zero-divisor graph of a commutative ring has also been studied in \cite{AFLL, SRR, SRR1, Rd1} and has been extended to non-commutative rings and semigroups in \cite{FL, Rd}.\\

The combinatorial properties of zero-divisors discovered in \cite{Bk} has also been studied in module theory. Recently in \cite{SR}, the elements of a module $M$ has been classified into \textit{full-annihilators}, \textit{semi annihilators} and \textit{star-annihilators}, see Definition 2.1 in section 2. For $M = R$, these elements are  the zero-divisors of a ring $R$, so the three simple graphs $ann_f(\Gamma(M))$, $ann_s(\Gamma(M))$ and $ann_t(\Gamma(M))$ corresponding to \textit{full-annihilators}, \textit{semi annihilators} and \textit{star-annihilators} in $M$ are natural generalizations of a zero-divisor graph introduced in \cite{AdLs}.\\

On the other hand, the study of essential ideals in a ring $R$ is a classical problem. For instance, Green and Van Wyk in \cite{GV} characterized essential ideals in certain class of commutative and non-commutative rings. The authors in \cite{A, KR} also studied essential ideals in $C(X)$ and topologically characterized the scole and essential ideals. Moreover, essential ideals also have been investigated in rings of measurable functions \cite{M} and $C^{*}$- algebras \cite{KP}. For more on essential ideals see \cite{HS, J, P}.\\

We call any subset of $M$ as an object. By combinatorial object we mean an object which can be put into one-to-one correspondence with a finite set of integers and by an algebraic object we mean a combinatorial object which is also an algebraic structure. The main objective of this paper is to study combinatorial objects, objects arising from modules and the graphs with vertex set as objects and combinatorial objects.\\

For an $R$-module $M$ and $x\in M$, set $[x : M] =\{r\in R ~|~ rM\subseteq Rx\}$, which clearly is an ideal of $R$ and an annihilator of the factor module $M/Rx$, where as the annihilator of $M$ is $[0 : M]$. In section 2, we study the correspondence of ideals in $R$, submodules of $M$ and the elements of an object $\widehat{A_{f}(M)}$, and characterize all finite abelian groups (Proposition 2.2). We further show (Theorem 2.11) that an $R$-module $[x : M]$ is injective if and only if $R$ is non-singular and the radical of $R/[x : M]$ is zero. In section 3, we examine two different equivalence relations on the elements of an object $\widehat{A_{f}(M)}$ and discuss (Theorem 3.6) the adjacencies of vertices in the graph $ann_f(\Gamma(M))$. Furthermore, we explore the equivalence relations (Theorem 3.8, Theorem 3.10) to establish the structure of a module $M$ and the annihilating graph $ann_f(\Gamma(M))$.
Finally, in section 4, we study some applications of annihilating graphs (Theorem 4.2) and extract certain module theoretic properties from these graphs. Further, we discuss the annihilating graphs arising from the tensor product (Theorem 4.7) and show that $ann_f(\Gamma(M \otimes_R T^{-1}R)) \cong ann_f (\Gamma(M))$ for every $R$-module $M$.\\

We conclude this section with some notations, which are mostly standard and will be used throughout this research article.\\

Throughout, $R$ is a commutative ring (with $1 \neq 0$) and all modules are unitary unless otherwise stated. A submodule $N$ is said to be an essential submodule of $M$ if it intersects non-trivially with every nonzero submodule of $M$. $[N : M] = \{r\in R~ |~ rM \subseteq N\}$ denotes an ideal of ring $R$. The symbols $\subseteq$  and $\subset$ has usual set theoretic meaning as containment and proper containment of sets. We will denote the ring of integers by $\mathbb{Z}$, positive integers by $\mathbb{N}$ and the ring of integers modulo $n$ by $\mathbb{Z}_n$. For basic definitions from graph theory we refer to [\cite{Dstl, Wst}], and for ring theory and module theory we refer to [\cite{AtDl, CE, Kp, W}].

\section{\bf Essential ideals determined by elements of an object ${\bf \widehat{A_{f}(M)}}$}
In this section, we discuss the correspondence of essential ideals in $R$, submodules of $M$ and the elements of an object $\widehat{A_{f}(M)}$. We characterize essential ideals corresponding to $\mathbb{Z}$-modules. We show that if $M$ is not simple $R$-module, then an ideal $[x : M]$, $x\in \widehat{A_{f}(M)}$ considered as an $R$-module is injective. We study essential ideals corresponding to the vertices of graph  $ann_f(\Gamma(M))$ over hereditary and regular rings.\\

We recall a definition concerning \textit{full-annihilators}, \textit{semi annihilators} and \textit{star-annihilators} of a module $M$.\\

\noindent{\bf Definition 2.1.} An element $x\in M$ is a,

(i) \emph{full-annihilator}, if either $x = 0$ or $[x : M][y : M]M = 0$, for some nonzero $y\in M$ with $[y : M] \neq R$,

(ii) \emph{semi-annihilator}, if either $x = 0$ or $[x : M] \neq 0$ and $[x : M][y : M]M = 0$, for some nonzero $y\in M$ with $ 0 \neq[y : M] \neq R$,

(iii) \emph{star-annihilator}, if either $x = 0$ or $ann(M) \subset  [x : M_R]$ and $[x : M][y : M]M = 0$, for some nonzero $y\in M$ with  $ann(M) \subset[y : M] \neq R$.\\

\noindent We denote by $A_f(M)$, $A_s(M)$ and $A_t(M)$ respectively the objects of full-annihilators, semi-annihilators and star-annihilators for any module $M$ over $R$. We set $\widehat{A_{f}(M)} = A_f(M)\backslash \{0\}$, $\widehat{A_{s}(M)} = A_s(M)\backslash \{0\}$ and $\widehat{A_t(M)} = A_t(M)\backslash \{0\}$.\\

In \cite{SR} authors introduced annihilating graphs arising from modules over commutative rings called as full-annihilating, semi-annihilating and star-annihilating graphs denoted by $ann_f(\Gamma(M))$, $ann_s(\Gamma(M))$ and $ann_t(\Gamma(M))$ respectively. The vertices of annihilating graphs are elements of objects $\widehat{A_{f}(M)}$, $\widehat{A_{s}(M)}$ and $\widehat{A_{t}(M)}$, and two vertices $x$ and $y$ are adjacent if and only if $[x : M][y : M]M = 0$. By Definition 2.1, we see that there is a correspondence of ideals in $R$, submodules of $M$ and the elements of objects $\widehat{A_{f}(M)}$, $\widehat{A_{s}(M)}$ and $\widehat{A_{t}(M)}$. Furthermore, the containment $ann_t(\Gamma(M)) \subseteq ann_s(\Gamma(M)) \subseteq ann_f(\Gamma(M))$ as induced subgraphs is clear, so our main emphasis is on the object $\widehat{A_{f}(M)}$ and the full-annihilating graph $ann_f(\Gamma(M))$. However, one can study these objects and graphs separately for any module $M$. Note that if $M$ is a finite module over $R$ or the graph $ann_f(\Gamma(M))$ is finite, then the objects $\widehat{A_{f}(M)}$, $\widehat{A_{s}(M)}$ are combinatorial with $|\widehat{A_{f}(M)}| = |\widehat{A_{s}(M)}|$ and the annihilating graphs $ann_f(\Gamma(M))$, $ann_s(\Gamma(M))$ coincide, where as the graph $ann_t(\Gamma(M))$ with vertex set as combinatorial object $\widehat{A_{t}(M)}$ may be different.\\

Let $G$ be any finite $\mathbb{Z}$-module. Clearly, $G$ is a finite abelian group. Below, we discuss the correspondence of ideals in $\mathbb{Z}$ and the elements of an object $\widehat{A_f(G)}$. We study cases of finite abelian groups where the essential ideals corresponding to the submodules generated by the vertices of graph $ann_f(\Gamma(M))$ are same and the submodules determined by these vertices are isomorphic.\\

The following is an interesting result in which a combinatorial object completely determines an algebraic object.\\

\noindent{\bf Proposition 2.2.} Let $G$ be a finite $\mathbb{Z}$-module. Then for each $x\in \widehat{A_f(G)}$, $[x : M]$ is an essential ideal if and only if $G$ is a finite abelian group without being simple.

\noindent{\bf Proof.} For all $x\in \widehat{A_f(G)}$, we have $[x : G] = n\mathbb{Z}$, $n\in\mathbb{N}$. It is clear that $n\mathbb{Z}$ intersects non-trivially with any ideal $m\mathbb{Z}$, $m\in\mathbb{N}$ in $\mathbb{Z}$.

For the converse, observe that among all finite abelian groups $\widehat{A_f(G)} = \emptyset$ if and only if $G$ is simple.\qed\\

Recall that a graph $\Gamma$ is said to be a complete if there is an edge between every pair of distinct vertices. A complete graph with $n$ vertices is denoted by $K_n$.\\

\noindent{\bf Remark 2.3.} For a finite abelian group $\mathbb{Z}_p\oplus \mathbb{Z}_p$, where $p \geq 2$ is prime, the essential ideals $[x : M]$, $x\in \widehat{A_f(\mathbb{Z}_p\oplus \mathbb{Z}_p)}$ corresponding to the submodules of $\mathbb{Z}_p\oplus \mathbb{Z}_p$ generated by elements of $\widehat{A_f(\mathbb{Z}_p\oplus \mathbb{Z}_p)}$ are same. In fact $[x : M] = ann(\mathbb{Z}_p\oplus \mathbb{Z}_p)$ for all $x\in \widehat{A_f(\mathbb{Z}_p\oplus \mathbb{Z}_p)}$. Furthermore, the abelian group $\mathbb{Z}_p\oplus \mathbb{Z}_p$ is a vector space over field $\mathbb{Z}_p$ and all one dimensional subspaces are isomorphic. So, the submodules generated by elements of $\widehat{A_f(\mathbb{Z}_p\oplus \mathbb{Z}_p)}$ are all isomorphic. For a finite abelian group $\mathbb{Z}_p\oplus \mathbb{Z}_q$, where $p$ and $q$ are any two prime numbers, the essential ideals determined by each $x\in \widehat{A_f(\mathbb{Z}_p\oplus \mathbb{Z}_q)}$ are either $p\mathbb{Z}$ or $q\mathbb{Z}$.\\

\noindent{\bf Remark 2.4.} Proposition 2.2 is not true for all $\mathbb{Z}$-modules. Consider a $\mathbb{Z}$-module $M = \mathbb{Z}\oplus \mathbb{Z}\oplus \cdots \oplus \mathbb{Z}$, which is a direct sum of $n$ copies of $\mathbb{Z}$. It is easy to verify that $\widehat{A_f(M)} = \widehat{M}$ with $[x : M][y : M]M = 0 $ for all $x, y\in M$, which implies $ann_f(\Gamma(M))$ is a complete graph. The cyclic submodules generated by the vertices of $ann_f(\Gamma(M))$ are simply the lines with integral coordinates passing through the origin in the hyper plane $\mathbb{R}\oplus \mathbb{R}\oplus \cdots \oplus \mathbb{R}$ and these lines intersect at the origin only. It follows that for each $x \in M$, $[x : M]$ is not an essential ideal in $\mathbb{Z}$, in fact $[x : M]$ is a zero-ideal in $\mathbb{Z}$.\\

Using the description given in Remark 2.4, we now characterize all the essential ideals corresponding to $\mathbb{Z}$-modules determined by  elements of $\widehat{A_f(M)}$.\\

\noindent{\bf Proposition 2.5.} If $M$ is any $\mathbb{Z}$-module, then $[x : M]$ is an essential ideal if and only if $[x : M]$ is non-zero for all $x\in \widehat{A_f(M)}$.

\noindent{\bf Proof.} $[x : M]$ is an ideal in $\mathbb{Z}$ for each $x\in \widehat{A_f(M)}$.\qed\\

For any $R$-module $M$, it would be interesting to characterize essential ideals $[x : M]$, $x\in \widehat{A_f(M)}$ corresponding to submodules determined by elements of $\widehat{A_f(M)}$ (or vertices of the graph $ann_f(\Gamma(M))$) such that the intersection of all essential ideals is again an essential ideal. It is easy to see that a finite intersection of essential ideals in any commutative ring is an essential ideal. But an infinite intersection of essential ideals need not to be an essential ideal, even a countable intersection of essential ideals in general is not an essential ideal as can be seen in \cite{A}. If $ann_f(\Gamma(M))$ is a finite graph, then by [Theorem 3.3, \cite{SR}] $M$ is finite over $R$, so the submodules determined by the vertices of graph are finite and therefore the ideals corresponding to submodules are finite in number and we conclude the intersection of essential ideals $[x : M]$, $x\in \widehat{A_f(M)}$ in $R$ is an essential ideal. Motivated by \cite{A}, we have have the following question regarding essential ideals corresponding to submodules $M$ determined by vertices of the graph $ann_f(\Gamma(M))$.\\

\noindent{\bf Question 2.6.} Let $M$ be an $R$-module. For $x\in \widehat{A_f(M)}$, characterize essential ideals $[x : M]$ in $R$ such that their intersection is an essential ideal.\\ 

 The Question 3.13 is true if every submodule of $M$ is cyclic with nonzero intersection.\\
  
\noindent{\bf Proposition 2.7.} Let $M$ be an $R$-module such that every submodule of $M$ is cyclic over $R$. For $x\in \widehat{A_f(M)}$, if the submodule generated by $x$ intersects non-trivially with every other submodule of $M$, then $[x : M]$ is an essential ideal in $R$. 

\noindent{\bf Proof.} Assume $\bigcap\limits_{x\in M} Rx \neq 0$. If $\widehat{A_f(M)} = \phi$, then $M$ is simple, a contradiction. Let $x\in \widehat{A_f(M)}$ and let $Rx$ be the submodule generated by $x$. Since $Rx$ intersects non-trivially with every other submodule, so there exist $y\in \widehat{A_f(M)}$ such that $Rx \cap Ry \neq 0$. It suffices to prove the result for $Rx \cap Ry$. Let $z\in Rx \cap Ry$ and let $[x : M]$, $[y: M]$,    $[z : M]$ be ideals of $R$ corresponding to submodules $Rx$, $Ry$ and $Rz$. Then $[z : M]\subseteq [x : M]\cap [y: M] \neq 0$, which implies $[x : M]$ intersects non-trivially with every nonzero ideal corresponding to the submodule generated by an element of $\widehat{A_f(M)}$. For any other ideal $I$ of $R$, it is clear that $IM = \{\sum\limits_{finite}am : a\in I, ~m\in M\} = Ra$ for some $a\in M$. Thus $I$ corresponds to the cyclic submodule generated by $a\in M$. It follows that $[x : M]\cap I \neq 0$, for every nonzero ideal of $R$ and we conclude that $[x : M]$ is an essential ideal for each $x\in \widehat{A_f(M)}$.\qed\\

The converse of  Proposition 2.7 is not true in general. We can easily construct examples from $\mathbb{Z}$-modules such that an ideal corresponding to the submodule generated by some element of an object $\widehat{A_f(M)}$ is an essential ideal, but the intersection of all submodules determined by elements of $\widehat{A_f(M)}$ is empty. However, if every ideal $[x : M]$, $x\in \widehat{A_f(M)}$ corresponds to an essential submodule of $M$, then we have the nonzero intersection.\\ 

\noindent{\bf Corollary 2.8.} Let $M$ be an $R$-module. For $x\in \widehat{A_f(M)}$, if the cyclic submodule $Rx$ intersects with every other cyclic submodule of $M$ non-trivially, then $[x : M]$ is an essential ideal in $R$.\\

\noindent{\bf Corollary 2.9.} The intersection $\bigcap\limits_{x\in \widehat{A_f(M)}}[x : M]$ is an essential ideal in $R$ if and only if every submodule of $M$ is essentially cyclic over $R$.\\

For an $R$-module $M$, let $Z(M)$ denote the following,
\begin{center}
$Z(M) = \{m\in M : ~ann(m) ~is ~an ~essential ~ideal ~in ~R\}$.
\end{center}
If $Z(M) = M$, then $M$ is said to be singular and if $Z(M) = 0$, then $M$ is said to be non-singular. By $rad(M)$, we denote the intersection of all maximal submodules of $M$. So, $rad(R)$ is the Jacobson radical of a ring $R$. The socle of an $R$-module $M$ denoted by $Soc(M)$ is the sum of simple submodules or equivalently the intersection of all essential submodules. To say that $Soc(M)$ is an essential socle is equivalent to saying that every cyclic submodule of $M$ contains a simple submodule of $M$. An essential socle of $M$ is denoted by $essoc(M)$.\\

In the following result, we consider singular simple $R$-modules (ideals) which are injective, and obtain some properties of essential ideals corresponding to the submodules generated by elements of $\widehat{A_f(M)}$.\\

\noindent{\bf Lemma 2.10.} For $x\in \widehat{A_f(M)}$, if $[x : M]$ is an essential ideal, then $R/[x : M]$ is a singular $R$-module.

\noindent{\bf Proof.} The proof simply follows by using definition of a singular module.\qed\\

\noindent{\bf Theorem 2.11.} Let $M$ be an $R$-module with $essoc(M)\neq 0$ and $\bigcap\limits_{x\in M} Rx \neq 0$. Then every singular simple $R$-module $[x : M]$, $x\in \widehat{A_f(M)}$ is injective if and only if $Z(R) = 0$ and $rad(R/[x : M]) = 0$.

\noindent{\bf Proof.} We have $essoc(M)\neq 0$ and $\bigcap\limits_{x\in M} Rx \neq 0$, so $\widehat{A_f(M)} \neq 0$. Therefore corresponding to every cyclic submodule generated by elements of $\widehat{A_f(M)}$ we have an ideal in $R$. Suppose all singular simple $R$-modules $[x : M]$, $x\in \widehat{A_f(M)}$ are injective. If for some $z\in \widehat{A_f(M)}$, $I = [z : M]\subseteq Z(R)$ is a simple $R$-module, then $Z(I) = I$. This implies $I$ is injective and thus a direct summand of $R$. However, the set $Z(R)$ is free from nonzero idempotent elements. Therefore $I = 0$ and it follows that  $Z(R) = 0$. For $x\in \widehat{A_f(M)}$, $A = [x : M]$ is an essential ideal of $R$. Thus by Lemma 2.10, $R/A$ is a singular module and so is every submodule of $R/A$. Therefore every simple submodule of $R/A$ is injective, which implies that every simple submodule is excluded by some maximal submodule. Thus we conclude that $rad(R/A) = 0$.

For the converse, we again consider the correspondence of cyclic submodules of $M$ and ideals of $R$. Let $\tilde{I}$ be a singular simple $R$-module corresponding to the submodule of $M$. In order to show that $\tilde{I}$ is injective, we must show that for every essential ideal $A$ in $R$ corresponding to the submodule determined by an element $x\in\widehat{A_f(M)}$, every $\varphi\in Hom_R(A, \tilde{I})$ has a lift $\psi\in Hom_R(R, \tilde{I})$ such that the following diagram commutes
\begin{center}
\begin{tikzcd}
 A \arrow{r}{f} \arrow[swap]{dr}{\varphi} & R \arrow{d}{\psi} \\
     & \tilde{I}
\end{tikzcd}
\end{center}
Let $K = ker(\varphi)$. We claim, $K$ is an essential ideal of $R$, for if $K \cap J = 0$ for some nonzero ideal $J$ of $R$. Then,
\begin{equation*}
I^{*} = J \cap A \neq 0 ~and ~I^{*} \cap K = 0.
\end{equation*}
This implies $I^{*}\subseteq \varphi(I^{*}) \subseteq \tilde{I}$, a contradiction, since $\tilde{I}$ is a singular simple submodule and $Z(R) = 0$. It is clear that if $\mu \neq 0$, $\varphi$ induces an isomorphism $\mu : A/K \longrightarrow \tilde{I}$. So, $A/K$ is a simple $R$-submodule of $R/K$. By our assumption $rad(R/K) = 0$, so there is a maximal submodule $M/K$ such that $R/K = A/K \oplus M/K$. Let $g : R \longrightarrow R/K$ be a canonical map and let $p : R/K \longrightarrow A/K$ be a projection map. Then we have $pg : R \longrightarrow A/K$. Therefore the composition $h = \mu pg$ is the required lift such that the above diagram commutes.\qed\\

In the rest of this section, we discuss some interesting consequences of preceding theorem.

\noindent{\bf Theorem 2.12.} Let $M$ be an $R$-module with $essoc(M)\neq 0$, $\bigcap\limits_{x\in M} Rx \neq 0$ and every singular simple $R$-module $[x : M]$, $x\in \widehat{A_f(M)}$ is injective. Then every ideal $[x : M]$ is an intersection of maximal ideals, $J(R)^{2} = 0$ and $[x : M]^{2} = [x : M]$. 

\noindent{\bf Proof.} For any $x\in \widehat{A_f(M)}$, $[x : M]$ is an essential ideal in $R$. Therefore, $J(R) \subseteq [x : M]$, since $J(R)$ is contained in every essential ideal of $R$. On the other hand intersection of all essential ideals in $R$ is scole of $R$, therefore $J(R) \subseteq Soc(R)$. This implies $J(R)^{2} = 0$ and $[x : M]$ is an intersection of maximal ideals in $R$. Suppose that $[x : M]^{2} \neq [x : M]$, for an essential ideal $[x : M]$ of $R$. By Theorem 2.11, $Z(R) = 0$ and therefore for every essential ideal $I$, $I \subseteq I^2$. In particular, $[x : M] \subseteq [x : M]^2$ for each $x\in \widehat{A_f(M)}$. It follows that $[x : M]^2$ is an essential ideal and is the intersection of maximal ideals in $R$. Finally, if $y\in [x : M]^2$, $y\notin [x : M]$, there is some maximal ideal $P$ of $R$ such that $[x : M] \subseteq P$, $y\notin P$. Then $R = Ry + P$, that is, $1 = ry + m$. This implies $y = yry + ym\in P$, a contradiction and hence we conclude that $[x : M]^{2} = [x : M]$.\qed\\

\noindent{\bf Corollary 2.13.} Let $M$ be an $R$-module, where $R$ is hereditary. For $x\in \widehat{A_f(M)}$, if $[x : M]$ is an essential ideal of $R$ and $J(R)^2 = 0$, then every singular simple $R$-module $[x : M]$ is injective.

\noindent{\bf Proof.} Let $R$ be hereditary. Then from \cite{CE}, the exact sequence,
\begin{equation*}
 0 \longrightarrow ann(x) \longrightarrow R \longrightarrow Rx \longrightarrow 0
\end{equation*}
splits for any $x\in R$. Since $J(R)^{2} = 0$ and $R/J(R)$ is an artinian ring. Therefore, $J(R) \subseteq Soc(R)$. But any essential ideal of $R$  contains $Soc(R)$. So, $J(R) \subseteq [x : M]$. This implies $R/[x : M]$ is completely reducible $R$-module and therefore $rad(R/[x : M]) = 0$. Thus by Theorem 2.11, every singular simple $R$-module $[x : M]$ is injective.\qed\\

Next, we consider modules over regular rings.\\

\noindent{\bf Theorem 2.14.} Let $M$ be an $R$-module such that every submodule of $M$ is cyclic over $R$ and $\bigcap\limits_{x\in M} Rx \neq 0$. Then the following are equivalent.

(i) $R$ is regular

(ii) $A^2 = A$ for each ideal $A$ of $R$

(iii) $[x : M]^{2} = [x : M]$ for each $x\in \widehat{A_f(M)}$

\noindent{\bf Proof.} The equivalence of (i) and (ii) is clear and certainly (ii) implies (iii). Thus we just need to show that (iii) implies (ii). By Theorem 2.11, $[x : M]$ is an essential ideal for each $x\in \widehat{A_f(M)}$. Suppose $[x : M]^{2} = [x : M]$. If $A$ is non-essential ideal of $A$, choose $J$ to be maximal ideal of $R$ such that $A\cap J = 0$, then $A + J$ is an essential ideal of $R$. Therefore again by Theorem 2.11, $A + J$ correspond to some submodule of $M$ and we have $A + J = [z : M]$ for some $z\in M$. So,
\begin{equation*}
(A + J)^2 = A^2 + J^2 = A + J.
\end{equation*}
If $x \in A$, then $x = \sum\limits_{finite}ab ~+  \sum\limits_{finite}mn$, where $a, b\in A$ and $m, n\in J$. Therefore,

\begin{equation*}
x - \sum\limits_{finite}ab = \sum\limits_{finite}mn \in A\cap J = 0.
\end{equation*}
This implies $x\in A^2$ and we conclude that $A = A^2$.\qed\\

\noindent{\bf Corollary 2.15.} Let $M$ be an $R$-module with $essoc(M)\neq 0$ and $\bigcap\limits_{x\in M} Rx \neq 0$. Then every singular simple $R$-module $[x : M]$, $x\in \widehat{A_f(M)}$ is injective if and only if $R$ is regular.

\noindent{\bf Proof.} By Theorem 2.12, if every singular simple $R$-module $[x : M]$ is injective, then for $x\in \widehat{A_f(M)}$, $[x : M]^{2} = [x : M]$. Therefore by Theorem 2.14, $R$ is regular. If $R$ is regular, then by [Theorem 6, \cite{RZ}] every singular simple $R$-module is injective.\qed\\

\section{\bf Equivalence relations on the elements of ${\bf\widehat{A_f(M)}}$, ${\bf\widehat{A_s(M)}}$, ${\bf\widehat{A_t(M)}}$}

This section is devoted for the study of equivalence relations defined on elements of objects $\widehat{A_f(M)}$, $\widehat{A_s(M)}$ and $\widehat{A_t(M)}$. We define two equivalence relations on $\widehat{A_f(M)}$ called as neighbourhood similar relation (combinatorial relation) and the submodule similar relation (algebraic relation). We investigate the equivalence of these two relations on $\widehat{A_f(M)}$ and study the conditions for any two elements of $\widehat{A_f(M)}$ to be adjacent in $ann_f(\Gamma(M))$. We explore the neighbourhood similar relation on $\widehat{A_f(M)}$ to establish the structure of module $M$ and the full-annihilating graph 
$ann_f(\Gamma(M))$.\\

\noindent{\bf Definition 3.1.} For an $R$-module $M$, two distinct elements $m_1, m_2\in \widehat{A_f(M)}$ are submodule similar $(\sim_{M})$ denoted by $m_1\sim_{M}m_2$ if and only if the submodules $ann(m_1)M$ and $ann(m_2)M$ of $M$ coincide, that is,
\begin{center}
$m_1\sim_{M}m_2$ if and only if $ann(m_1)M = ann(m_2)M$. 
\end{center}
Clearly, $\sim_{M}$ is an equivalence relation on $\widehat{A_f(M)}$ and the equivalence classes of any $m\in \widehat{A_f(M)}$ is denoted by,
\begin{center}
$[m]_M = \{m^{*}\in M : m^{*}\sim_{M}m\}$. 
\end{center}
Analogously, we can define the submodule relation on the elements of objects $\widehat{A_s(M)}$ and $\widehat{A_t(M)}$.\\ 

Let $\Gamma$ be any connected graph. A neighbour of any vertex $v$ in $\Gamma$ is a vertex adjacent to $v$. $N(v)$ denotes the set all neighbours of $v$ and $N[v] = N(v) \cup \{v\}$. The study of neighbourhoods of vertices in a connected graph $\Gamma$ is related to the symmetry of that graph. There is a close relationship, which is being discussed in the following definition between the neighbourhoods and the distance similar classes of vertices defined in \cite{SR1} .\\

\noindent{\bf Definition 3.2.} For a connected graph $\Gamma$, two distinct vertices $u, v\in V(\Gamma)$ on a vertex set $V(\Gamma)$ are neighbourhood similar $(\sim_{nbd})$ denoted by $u \sim_{nbd} v$ if and only if $N(u) = N(v)$.\\ 

It can be easily checked that $\sim_{nbd}$ is an equivalence relation on $V(\Gamma)$. The neighbourhood similar equivalence class of a vertex $v$ is,
\begin{center}
$[v]_{\Gamma} = \{w\in V(\Gamma) : N(v) = N(w)\}$.
\end{center}

Two distinct vertices $a, b\in V(\Gamma)$ which are not neighbourhood similar can be identified as the vertices for which $ab \in E(\Gamma)$ with $N(a) \neq N(b)$ or $ab \notin E(\Gamma)$ with $N(a) \neq N(b)$, where $E(\Gamma)$ denotes the edge set of $\Gamma$.\\

If $\Gamma$ is a finite connected graph, then the neighbourhood similar  relation on $V(\Gamma)$ is a distance similar relation ($\sim_{d}$) defined in \cite{SR1} with two vertices $u, v\in V(\Gamma)$ are distance similar denoted by $u\sim_{d} v$ if $d(u, w) = d(v, w)$ for all $w \in V(\Gamma)\setminus\{u, v\}$. Clearly, vertices $u$ and $v$ are distance similar if either $uv \notin E(\Gamma)$ and $N(u) = N(v)$
or $uv \in E(\Gamma)$ and $N[u] = N[v]$.\\

The following example illustrates neighbourhood similar relation on a connected graph $\Gamma$.\\

\noindent{\bf Example 3.3.} Consider the graph $\Gamma$ of order $n = 10$, shown in Figure 1. Two of the five neighbourhood similar equivalence classes are $V_1 = \{v_6, v_5\}$ and $V_2 = \{v_7, v_8, v_9, v_{10}\}$. Each of the three remaining classes $\{v_1\}$, $\{v_4\}$ and $\{v_6\}$ consists of a single vertex.

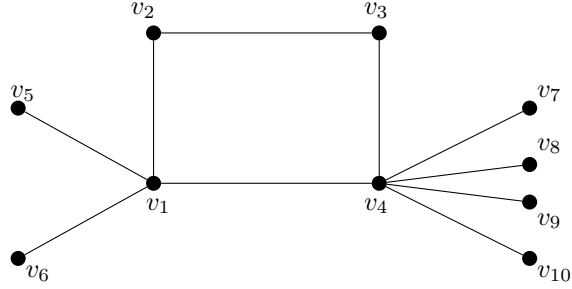
\begin{align*}
\begin{pgfpicture}{9cm}{-3cm}{6cm}{2cm}
\pgfnodecircle{Node1}[fill]{\pgfxy(7,0)}{0.1cm}
\pgfnodecircle{Node2}[fill]{\pgfxy(7,2)}{0.1cm}
\pgfnodecircle{Node3}[fill]{\pgfxy(10, 2)}{0.1cm}
\pgfnodecircle{Node4}[fill]{\pgfxy(10,0)}{0.1cm}
\pgfnodecircle{Node6}[fill]{\pgfxy(5.2, 1)}{0.1cm}
\pgfnodecircle{Node7}[fill]{\pgfxy(5.2,-1)}{0.1cm}
\pgfnodecircle{Node8}[fill]{\pgfxy(12,1)}{0.1cm}
\pgfnodecircle{Node9}[fill]{\pgfxy(12, .25)}{0.1cm}
\pgfnodecircle{Node10}[fill]{\pgfxy(12,-.25)}{0.1cm}
\pgfnodecircle{Node11}[fill]{\pgfxy(12,-1)}{0.1cm}
\pgfnodeconnline{Node1}{Node2}
\pgfnodeconnline{Node1}{Node4}
\pgfnodeconnline{Node1}{Node6}
\pgfnodeconnline{Node1}{Node7}
\pgfnodeconnline{Node2}{Node3}
\pgfnodeconnline{Node3}{Node4}
\pgfnodeconnline{Node4}{Node8}
\pgfnodeconnline{Node4}{Node9}
\pgfnodeconnline{Node4}{Node10}
\pgfnodeconnline{Node4}{Node11}
\pgfputat{\pgfxy(6.9,-.3)}{\pgfbox[left,center]{$v_1$}}
\pgfputat{\pgfxy(6.7, 2.3)}{\pgfbox[left,center]{$v_2$}}
\pgfputat{\pgfxy(9.8, 2.3)}{\pgfbox[left,center]{$v_3$}}
\pgfputat{\pgfxy(9.8,-.3)}{\pgfbox[left,center]{$v_4$}}
\pgfputat{\pgfxy(5.1, 1.2)}{\pgfbox[left,center]{$v_5$}}
\pgfputat{\pgfxy(5.3,-1.2)}{\pgfbox[left,center]{$v_6$}}
\pgfputat{\pgfxy(12.1, 1.2)}{\pgfbox[left,center]{$v_7$}}
\pgfputat{\pgfxy(12.1, .5)}{\pgfbox[left,center]{$v_8$}}
\pgfputat{\pgfxy(12.1, -.5)}{\pgfbox[left,center]{$v_9$}}
\pgfputat{\pgfxy(12.1, -1.2)}{\pgfbox[left,center]{$v_{10}$}}
\pgfputat{\pgfxy(3.4, -2)}{\pgfbox[left,center].{$Figure \hskip .1cm 1. \hskip .1cm Graph \hskip .1cm of \hskip .1cm order \hskip .1cm 10 \hskip .1cm with \hskip .1cm 5 \hskip .1cm neighbourhood \hskip .1cm similar \hskip .1cm classes$}}
\end{pgfpicture}
\end{align*}

Now, we explain the connection between neighbourhood similar and submodule similar equivalence relations on an object $\widehat{A_f(M)}$. In fact, we will see that the neighbourhood similar relation $\sim_{nbd}$, which is a combinatorial relation implies an algebraic relation $\sim_{M}$ and conversely. Moreover, we investigate the condition for any two elements of  $\widehat{A_f(M)}$ to be adjacent in $ann_f(\Gamma(M))$. We start with following two lemmas.\\

\noindent{\bf Lemma 3.4.} Let $M$ be an $R$-module. Then for $x, z, y\in \widehat{A_f(M)}$, $[x : M][z : M]M = 0$ if and only if $z \in ann(x)M$ or $z \in ann(y)M$ with $ann(x)M \bigcap ann(y)M = \{0\}$.

\noindent{\bf Proof.} The proof follows because of the fact that each element of an ideal $ann(y)$ in $R$ annihilate whole of $Ry$.\qed\\

\noindent{\bf Lemma 3.5.} If $M$ is an $R$-module and $x, y\in\widehat{A_f(M)}$ with $x \sim_{nbd} y$, then either $[x : M]^2M = [y : M]^2M = 0$ or both $[x : M]^2M \neq 0$ and $[y : M]^2M \neq 0$.

\noindent{\bf Proof.} Let $M$ be an $R$-module and let $x \sim_{nbd} y$ for $x, y\in\widehat{A_f(M)}$. Assume that $[x : M]^2M = 0$ and $[y : M]^2M \neq 0$. We consider the following cases.

Case 1. $x$ and $y$ are adjacent in $ann_f(\Gamma(M))$, that is $[x : M][y : M]M = 0$. Then it is easy to check that $[x+y : M][x : M]M = 0$, which implies $x+y \in N(x)$. Note here that $x+y \neq x$ and $x+y \neq y$. However, $x+y \notin N(y)$, since $[x+y : M][y : M]M \neq 0$. Therefore $N(x) \neq N(y)$, a contradiction.

Case 2. $x$ and $y$ are not adjacent in $ann_f(\Gamma(M))$, that is $[x : M][y : M]M \neq 0$. Then by case 1, we can easily find an element $z \in\widehat{A_f(M)}$ such that $[z : M][x : M]M = 0$, where as $[z : M][y : M]M \neq 0$. Therefore, $N(x)\neq N(y)$, again a contradiction.\qed\\

\noindent{\bf Theorem 3.6.} Let $M$ be an $R$-module with $|\widehat{A_f(M)}|$ being at least 3. Then for $x,y\in\widehat{A_f(M)}$ the following hold.

(i) If $M$ is a multiplication module, then $x \sim_{nbd} y$ if and only if  $x \sim_{M} y$.

(ii) If $x \sim_{nbd} y$, then $[x : M][y : M]M = 0$ if and only if $[x : M]^2M = [y : M]^2M = 0$.

\noindent{\bf Proof.} (i) Directly follows from Lemma 3.4 and the fact that every submodule of a multiplication module $M$ is of the form $IM$, for some ideal $I$ of $R$. 

(ii) Let $x \sim_{nbd} y$ and suppose $[x : M]^2M \neq 0$. Then by Lemma 3.5, $[y : M]^2M \neq 0$. Assume that $[x : M][y : M]M = 0$. There must be without loss of generality some $m\in N(x)$ with $m \neq x$, such that $m\in N(y)$. Then $[y : M][x+m : M]M = 0$, which implies either $x+m\in N(x)$ or $x+m = x$ or $x+m = y$. But, $x+m \neq x$ as $[x : M]^2M \neq 0$, also $x+m \notin N(x)$, since $[x : M][x+m : M]M \neq 0$ and $x+m \neq y$, since $y \in N(x)$. Lastly, if $x+m = 0$, then $[x : M]^2M = 0$. Thus, we have a contradiction and therefore, $[x : M][y : M]M \neq 0$.

Now, suppose $[x : M]^2M =  0$. Therefore by Lemma 3.5, $[y : M]^2M = 0$. Assume that $[x : M][y : M]M \neq 0$. Then by Lemma 3.4, $x \notin ann(x)M$ and $x \notin ann(y)M$. Since $|\widehat{A_f(M)}| \geq 3$, so there is some $z\in \widehat{A_f(M)}$ such that $[x : M][z : M]M = 0$. This implies $z = x$ or $z = y$, since $x \sim_{nbd} y$. But, $z \neq y$ as $ y \notin N(x)$. Therefore $z = x$ and we have for all $r\in [x : M]$, $x+ry \in ann(x)M$ or $x+ry \in ann(y)M$. However, $x+ry \notin ann(y)M$, since $[x + ry : M][y : M]M \neq 0$. On the other hand, $x+ry \in N(y)$, since $x \sim_{nbd} y$, a contradiction. If $x+ry = x$ or $x+ry = y$, then $[x : M][y : M]M = 0$,  which is again a contradiction.\qed\\

The preceding theorem is also true for the elements of objects $\widehat{A_s(M)}$ and $\widehat{A_t(M)}$, that is, the same adjacency relations hold for the vertices of semi-annihilating $ann_s(\Gamma(M))$ and the star-annihilating graph $ann_t(\Gamma(M))$.\\

Recall that an element $a\in R$ is said to be nilpotent if for some $n\in\mathbb{N}$, $a^n = 0$. We denote the set of all nilpotent elements of $R$ by $nil(R)$ called the nil radical of $R$ which is contained in every prime ideal of $R$. For a module $M$, we denote by $nil(M)$ the sub module which is contained in every prime submodule of $M$ that is,
\begin{center} $nil(M) = \bigcap\limits_{N\in Spec(M)}N = \{x\in M ~|~[x : M]x = 0\}$,\end{center}
where $N = \{x\in M ~|~ [x : M]x = 0\}$ and Spec(M) is the set of all prime submodules of $M$. Also, note that a  vertex $v\in V(\Gamma)$ is said to be a pendant vertex if there is only one vertex adjacent to it. That is, if degree of vertex a $v$, denoted by $deg(v)$ is 1.\\

In the remaining section, we study the neighbourhood classes of elements in $M$. Consider the elements of $M$ which are not neighbourhood similar in $ann_f(\Gamma(M))$, we investigate about the structure of a module $M$ if any two elements in $ann_f(\Gamma(M))$ are not neighbourhood similar. We make use of neighbourhood similar elements to study the nature of all annihilating graphs arising from $M$. Before, we discuss the results regarding neighbourhood variant and invariant classes we have the following key lemmas.\\

\noindent{\bf Lemma 3.7.} For any $R$-module $M$, $ann_f(\Gamma(M))$ is not an $n-gon$ for $n\geq 5$.

\noindent{\bf Proof.} Suppose $ann_f(\Gamma(M))$ is the graph with vertices $\{x_1, x_2, \cdots, x_n\}$ such that $x_1 - x_2, x_2 - x_3, \cdots, x_{n-2} - x_{n-1}, x_{n-1} - x_n,  x_n - x_1$ are the only adjacencies in $ann_f(\Gamma(M))$. Then, we have $[x_1 : M][x_2 : M]M = 0 = [x_1 : M][x_n : M]M$. This implies $[x_1 : M][x_2 + x_n : M]M = 0$. It follows that $x_2 + x_n$ is either $x_1, x_2, \cdots, x_{n-1}$ or $x_n$. A simple check yields that $x_2 + x_n = x_1$ is the only possibility. Similarly, $x_1 + x_{n-1} = x_n$. Therefore, $x_n = x_2 + x_n + x_{n-1}$. So, 
\begin{equation*} [x_1 : M][x_n : M]M = [x_1 : M][x_2 : M]M + [x_1 : M][x_n : M]M + [x_1 : M][x_{n-1} : M]M.
\end{equation*} 
This implies $[x_1 : M][x_{n-1} : M]M = 0$, a contradiction, and hence we conclude that $ann_f(\Gamma(M))$ is not an $n-gon$ for $n\geq 5$.\qed\\

\noindent{\bf Lemma 3.8.} Let $M$ be an $R$-module. Then for all $x\in M$, $[x : M]x = 0$ if and only if $[x : M]^{2}M = 0$.

\noindent{\bf Proof.} The proof follows because of the fact that an ideal $[x : M]$ of $R$ annihilate every multiple of $x$.\qed\\

\noindent{\bf Lemma 3.9.} For an $R$-module $M$ with $nil(M) \neq 0$, let $ann_f(\Gamma(M))$ be the full annihilating graph such that $ann_f(\Gamma(M))$ contains no cycle of length 3. Then the following hold.

(i) If for every pair $y, z\in \widehat{A_f(M)}$ with $y \sim_{nbd} z$ there exists an element $x\in  \widehat{A_f(M)}$ such that $[x : M][z : M]M = 0 = [x : M][y : M]M$, then either $8\leq |M|\leq 16$ or $|M|\geq 17$ and $nil(M) = \{0, x\}$.

(ii) If for each $x\in \widehat{A_f(M)}$, there exists a pair $z, y\in \widehat{A_f(M)}\setminus\{x\}$ such that $y \sim_{nbd} z$ and $[x : M][y : M]M = 0 = [x : M][z : M]M$, then any $a\in \widehat{A_f(M)}$ with $[x : M][a : M]M = 0$ is a pendant vertex in $ann_f(\Gamma(M))$.

\noindent{\bf Proof.} (i) Let $nil(M) \neq 0$. We first show that for all $r\in[x : M]$, $r^nx = 0$ for some $n\in \mathbb{N} $, where $x\in nil(M)$. Consider a set $S = \{r^nx : n\in \mathbb{N}\}$, we must show that $0 \in S$. Suppose that $0 \notin S$ and let $N$ be a submodule of $M$ such that $N \cap S = \emptyset$. Then by Zorn's Lemma the collection $\Sigma = \{N : N\cap S = \emptyset\}$ contains a maximal element. Let $K$ be a maximal member of $\Sigma$. We claim that $[N : M]$ is a prime ideal of R. Clearly, $[N : M] \subset R$. Let $r_1, r_2 \in R$ and suppose $r_1r_2\in [N : M]$ with 
$r_1, r_2 \notin [N : M]$. Then $(r_1M + K) \notin \Sigma$ and $(r_2M + K) \notin \Sigma$. So $r^{n}_{1}x \in S \cap r_1M + K$ and $r^{n}_{2}x \in S \cap  r_2M + K$ for some $n_1, n_2 \in \mathbb{N}$. Therefore $r^{n_{1}+n_{2}}x \in K \cap S$, a contradiction. Thus $[N : M]$ is a prime ideal and by [Corollary 2.11, \cite{ES}], $K$ is a prime submodule of $M$. Therefore we have $rx\in N \cap S$, since $x\in [N : M]$, which is a contradiction and consequently $r^nx = 0$. By well ordering principle choose $n$ to be smallest such that $r^nx = 0$. Then for $n \geq 1$, $r^{n-1}x \neq 0$. 

Claim: $n \leq 3$. Suppose to the contrary that $n > 3$. Clearly, $r^{n-1}x, ~rx\in \widehat{A_f(M)}$, since $[r^{n-1}x : M][rx : M]M = 0$. So, there exist a vertex say $y\in \widehat{A_f(M)}$ such that $[r^{n-1}x : M][y : M]M = 0 = [y : M][rx : M]M$, which implies $[r^{n-1}x : M][y : M]M = 0 = [r^{n-1}x : M][rx : M]M$. Therefore $r^{n-1}x = y$ is the only possibility, since $r^{n-1}x \sim_{nbd}rx$. Similarly for each $i$, $1 \leq i \leq n-2$, we have $[r^{i}x : M][r^{n-1}x : M]M = 0$. For $m = r^{n-2}x + r^{n-1}x$, we see that $ann_f(\Gamma(M))$ contains a cycle of length 3, since $[r^{n-1}x : M][m : M]M = 0 = [r^{n-2}x : M][m : M]M = 0$ with $m \notin \{0, ~r^{n-1}x, ~r^{n-2}x\}$. Thus $n \leq 3$.

We consider the following cases for $n\leq 3$.

Case 1. $n = 3$. We show that $ann(r^2x)M$ is the unique maximal submodule of $M$ and $|M| = 8$ or $|M| = 16$. We first show that $|M| = 8$ or $|M| = 16$. By our claim above, $[r^{2}x : M][rx : M]M = 0$. If $0 \neq z\in ann(x)M$, then $[z : M][x : M]M = 0$. This implies $ann_f(\Gamma(M))$ contains a cycle of length 3 on vertices $z$, $rx$ and $r^2x$. Therefore $z = r^2x$ and $ann(x)M\subseteq \{0, r^2x\}$. In fact, $ann(x)M \subseteq Rr^{2}x$, since for all $s\in R$, $[sr^{2}x : M][rx : M]M = 0 = [sr^{2}x : M][r^2x : M]M$. Therefore $sr^{2}x  \in \{0, ~r^{n-1}x, ~r^{n-2}x\}$. If $sr^{2}x = rx$, then $r^{2}x = 0$, a contradiction. Thus $Rr^{2}x = \{0, r^{2}x\}$.

Moreover, $ann(r^{2}x)M \subseteq \{0, x, rx, r^{2}x, x + rx, x + r^{2}x, rx +        r^{2}x, x + rx + r^{2}x\}$. If $ z \in ann(r^{2}x)M$, then $r^{2}z \in ann(x)M\subseteq \{0, r^2x\}$. So, either $r^{2}z = 0$ or $r^{2}z = r^{2}x$. Therefore $[rz : M][rx : M]M = 0 = [rz : M][r^2x : M]M$ or $[rz - rx : M][rx : M]M = 0 = [rz - rx : M][r^2x : M]M$. This implies $rz\in \{0, rx, r^2x\}$ or $rz - rx\in \{0, rx, r^2x\}$. Let $r^{2}z = 0$. So, $rz \neq rx$ and therefore either $rz = 0$ or $r(z - rx) = 0$, which implies $[z : M][rx : M]M = 0 = [z : M][r^2x : M]M$ or $[z - rz : M][rx : M]M = 0 = [z - rz : M][r^2x : M]M$. Thus $z\in \{0, rx, r^{2}x, rx + r^{2}x \}$. Therefore we may assume that $r^2z = r^2x$, which implies $rz - rx \neq rx$. But, $rz - rx \in \{0, rx, r^{2}x\}$. So, either $rz - rx = 0$ or $rz - rx = r^2x$ and by similar argument as above $z \in \{x, r^{2}x, x + rx, x + rx + r^{2}x\}$.

\noindent If $[x : M][r^2x : M]M = 0$, then 
\begin{equation}
ann(r^2x)M = \{0, x, rx, r^{2}x, x + rx, x + r^{2}x, rx + r^{2}x, x + rx + r^{2}x\}. 
\end{equation}    
If $[x : M][r^2x : M]M \neq 0$, then
\begin{equation}
ann(r^2x)M = \{0, rx, r^{2}x, rx + r^{2}x\}   
\end{equation}
For (1), $|M| = 16$ and for (2), $|M| = 8$, since $[x : M]r^{2}M \neq 0$, so there are $t\in [x : M]$ and $m\in M$ such that $r^2tm \neq 0$. It is clear that $r^2tm = r^{2}x$. Let $m^{*}\in M$. Then $r^2tm^{*}\in Rr^{2}x = \{0, ~r^{2}x\}$. If $r^2tm^{*} = 0$, then $m^{*}\in ann(r^2x)M$ and if $r^2tm^{*} =      r^{2}x$, then $m^{*} - m\in ann(r^2x)M$.

Let $ K = ann(r^{2}x)M$. Clearly, $K \subset M$ and $Rr^{2}x \cong R/ann(r^{2}x)$. Therefore $ann(r^{2}x)$ is a maximal ideal, since $Rr^{2}x = \{0, r^{2}x\}$. Thus it follows by [Theorem 2.5, \cite{ES}] that $K$ is a maximal submodule. Further, $K \subseteq Rx \subseteq nil(M) \subseteq K$. Therefore $K = nil(M)$ is the unique maximal submodule of $M$. If $A_f(M) \subseteq K$, then $A_f(M) = K$, which implies $ann_f(\Gamma(M))$ is a star graph.  

Case 2. $n = 2$. We show that $|M| \leq 12$. If $[x : M]^{2}x = [x : M]^{3}M \neq 0$, then there exist two elements $r, s\in [x : M]$ such that $rsx \neq 0$. Further, there are $m\in M$ and $t\in [x : M]$ such that $rstm \neq 0$. However, $r^{2}x = s^{2}x = t^{2}x = 0$ and there is some $a\in\widehat{A_f(M)}$ such that $[a : M][rx : M]M = 0$. It is easy to check that $Rrx \subseteq \{0, rx, a\}$ and $a = rsx$. This implies that $[srx : M][rx : M]M = 0$ and therefore $Rrx = \{0, ~rx, ~rsx\}$ and $ann(rx)M = \{0, ~rx, ~rsx\}$. Clearly, $[rstm : M][s_1x : M]M = 0$, $[rstm : M][s_1s_2x : M]M = 0$ with $rstm \neq s_1x$. Thus $rstm = rsx$ and it is clear that $[rsm: M][rx : M]M = 0$ and $[rsm : M][rsx : M]M = 0$. But, $rsm \neq s_1x$ and $rsm \neq rsx$, a contradiction. Therefore $[x : M]^{2}x = 0$.

Moreover, we have $rsm \neq 0$, for some $s\in[x : M]$ and $m\in M$, since $r^{2}x = 0$, $rx \neq 0$. For $x, y\in \widehat{A_f(M)}$ with $[x : M][y : M]M = 0$, there is $rx\in \widehat{A_f(M)}$ such that $[x : M][rx : M]M = 0$ and $[rx : M][y : M]M = 0$. So, $rx = x$ or $rx = y$. If $rx = x$, then $rx = 0$, which is not possible. Therefore $rx = y$ and we have $[x : M][rx : M]M = 0$. Let $z\in ann(x)M$. Then $z\in \{0, x, rx\}$, since $[x : M][rx : M]M = 0$. If $z = x$, then $x[x : M] = [x : M]^{2}M = 0$, which is a contradiction. Thus $ann(x)M = \{0, rx\}$ and in fact $Rrx = \{0, rx\}$. On the other hand, $rm\in \widehat{A_f(M)}$, so there exist $b\in \widehat{A_f(M)}$ such that $[rm : M][b : M]M = 0$. But, $[rsm : M][rm : M]M = 0$ and $[rsm : M][b : M]M = 0$, therefore $rsm = b$ is the only possibility, so $[rsm : M][rm : M]M = 0$. Also, $[rsm : M][rx : M]M = 0$ and $[rsm : M][x : M]M = 0$. Thus by the same reasoning as above we have $rsm = rx$. Let $t\in ann(rx)M$. Therefore $rt\in  ann(x)M = \{0, rx\}$. If $rt = 0$, then $[t : M][rsm : M]M = 0 = [t : M][rm : M]M$. If $rt = rx$, then $[t - x : M][rsm : M]M = 0 = [t - x : M][rsm : M]M$.  Therefore $ann(rx)M = \{0, rm, rx, x + rm, x + rx\}$. By similar argument as in case 1, it can be shown that $|M| \leq 12$ and $ann(rx)M = nil(M)$ is the unique maximal submodule of $M$.
 
Case 3. $n = 1$. If $[x : M]x = [x : M]^{2}M \neq 0$, then by cases 1 and 2 we have $8\leq |M|\leq 16$. Assume that $[x : M]x = 0$, we show that either $|M| = 9$ or $nil(M) = \{0, x\}$ with $2x = 0$. Let $x\in [x : M]M$. Then $x = \sum\limits_{i \in \Lambda}r_im_i$, where $\Lambda$ is finite, $r_i\in [x : M]$ and $m_i\in M$ with $1\leq i\leq |\Lambda|$. For $x\in \widehat{A_f(M)}$, there is $rx\in \widehat{A_f(M)}$ such that $[x : M][y : M]M = 0$. So, $Rx \subseteq \{0, ~x, ~y\}$. If $x \neq r_im_i$ for all $i$, $1\leq i\leq |\Lambda|$, then $r_im_i\in Rx$. Therefore $r_im_i = y$ for all $i$, $1\leq i\leq |\Lambda|$. Thus it follows that $x = rm$ for some $r\in [x : M]$ and $m\in M$ with $rm \neq 0$. Clearly, $x + x\in Rx\subseteq \{0, ~x, ~y\}$. If $x + x \neq 0$, then $Rx = \{0, x, 2x\}$, $[x : M][2x : M]M = 0$ and $ann(x)M = \{0, x, 2x\}$. Thus for all $m^{*}\in M$, $rm^{*}\in Rx$ and we have, 
\begin{center}$[m^{*} : M][x : M]M = 0 = [m^{*} : M][2x : M]M$, \end{center} or 
\begin{center}$[m^{*}-m : M][x : M]M = 0 = [m^{*}-m : M][2x : M]M$,\end{center} or
\begin{center}$[m^{*}-2m : M][x : M]M = 0 = [m^{*}-2m : M][2x : M]M$.\end{center}  
By a similar argument as in case 1 it can be shown that $|M| = 9$ and $ann(x)M$ is the unique maximal submodule of $M$. Let $|M| \neq 9$. Then by the above argument we must have $2x = 0$. We show that $nil(M) = \{0, x\}$. If $0\neq z \in nil(M)$ with $z \neq x$, then $z[z : M] = [z : M]^{2}M = 0$. Therefore $z = sm_1$ for some $s\in [z : M]$ and $m_1\in M$ with $s^{2}m_1 \neq 0$. It is clear that $[x : M][x_1 : M]M = 0$ and $[z : M][z_1 : M]M = 0$ with $x_1, z_1 \in \widehat{A_f(M)}$, since $x, z \in \widehat{A_f(M)}$. Thus $Rx \subseteq \{0, x, x_1\}$ and $Rz \subseteq \{0, z, z_1\}$. If $0 \neq rsm\in Rx$ and $rsm\in Rz$, say $rsm = x\in Rz$, then $x = z_1$ which implies $[x : M][z : M]M = 0$. Therefore $rsm = 0$, a contradiction. If $rsm = x_1$, then $Rx =\{0, x, rsm\} = ann(x)M$ and by similar argument as above we see that $|M| = 9$, again a contradiction. Thus $rsm = 0$ and similarly $rsm_1 = 0$. Clearly, $x + z \neq x$, $x + z \neq z$ and $x + z \in\widehat{A_f(M)}$, so there is $t\in\widehat{A_f(M)}$ such that $[x + t : M][t : M]M = 0$. Also, $rt \in Rx \subseteq \{0, x, x_1\}$. If $rt = 0$, then $[x : M][t : M]M = 0 = [x : M][x + z : M]M$. That is, we have a cycle on vertices $x$, $t$ and $x + z$, a contradiction. If $rt = x_1$, then $Rx =\{0, x, rt\} = ann(x)M$ and by the same argument as above we have $|M| = 9$, again a contradiction. Therefore we may assume that $x = rt$ and $z = st$. Then
\begin{center}$[x + z : M][x : M]M = 0 = [x + z : M][t : M]M$.\end{center}
This implies $x + z = 0$. Thus we have a contradiction in every possible case and hence we conclude that $nil(M) = \{0, x\}$ with $2x = 0$.

(ii) By (i), we see that $nil(M) = \{0, x\}$ for some $0\neq x\in M$ with $2x = 0$, $[x : M]^{2}M = 0$ and $|M| \geq 17$. Let $x\in \widehat{A_f(M)}$. Then, there is $y\in \widehat{A_f(M)}$ such that $[x : M][y : M]M = 0$. Clearly, there is another $x + y\in \widehat{A_f(M)}$ with $x + y \neq x$ and $x + y \neq y$ such that $[x + y : M][y : M]M = 0$. Therefore, $x + y \sim_{nbd} y$, that is $N(x + y) = N(y)$. If there is some $t\in \widehat{A_f(M)}$ such that  $[x + y : M][t : M]M = 0$, then $[y : M][t : M]M = 0$. This implies $[x : M][t : M]M = 0$, a contradiction, since $ann_f(\Gamma(M))$ contains no cycle of length 3. Thus, $x + y = t$ or $y = t$. But, $ y\neq t$, otherwise $x = 0$, a contradiction, since $nil(M) = \{0, x\}$ with $0\neq x$. Therefore $x + y = t$ and consequently $t\in \widehat{A_f(M)}$ is a pendant vertex in $ann_f(\Gamma(M))$.\qed\\

Using the neighbourhood similar relation on the elements of objects  $\widehat{A_f(M)}$, $\widehat{A_s(M)}$ and $\widehat{A_t(M)}$, we now reveal the structure of a module $M$.\\

\noindent{\bf Theorem 3.8.} Let $M$ be an $R$-module and let $ann_f(\Gamma(M))$ be the full annihilating graph of $M$ which contains no cycle of length 3. If for each $x\in \widehat{A_f(M)}$, there exists a pair $z, y\in \widehat{A_f(M)}\setminus\{x\}$ such that $y \not\sim_{nbd} z$ and $[x : M][y : M]M = 0 = [x : M][z : M]M$, then $M = M_1 \oplus M_2$, for some submodules $M_1$ and $M_2$ of $M$.

\noindent{\bf Proof.} Let $x\in \widehat{A_f(M)}$. There are $z, y\in \widehat{A_f(M)}\setminus\{x\}$ such that $y \not\sim_{nbd} z$ and $[x : M][y : M]M = 0 = [x : M][z : M]M$. So, there is some $t\in \widehat{A_f(M)}$ such that $[t : M][y : M]M = 0$, where as $[t : M][z : M]M \neq 0$. Clearly, for all $r\in [z : M]$, $0\neq rt\in A_f(M)$ and we have $[rt : M][y : M]M = 0$. This implies $rt = x$ or $rt = y$, since $[x : M][y : M]M = 0$. Therefore $[y : M]^2M = 0 = y[y : M]$ or $[x : M]^2M = 0 = x[x : M]$. That is, $x\in nil(M)$ or $y\in nil(M)$. By Lemma 3.7, $nil(M) = \{0, m\}$ for some $0\neq m\in M$. Suppose that $x = m$. Then $[m : M][y : M]M = 0$ and $[m : M][t : M]M \neq 0$, otherwise $ann_f(\Gamma(M))$ contains a cycle of length three on $t$, $m$ and $y$. So, there is some $r_1\in [t : M]$ such that $r_1m \neq 0$. It is clear that $[r_1m : M]^2M = 0$, since $[m : M]^2M = 0$. Therefore for all $r_1\in  [m : M]$, $r_1m = m$. Let $s = r_1t - t$. Then $[s : M][y : M]M = 0$, since  $[t : M][y : M]M = 0$. For $r_1\in  [m : M]$, we see that $us = ur_1t - ut = ut - t = 0$. Therefore $[m : M][s : M]M = 0$. This implies $s\in \{0, y, m\}$, since $[m : M][y : M]M = 0$. We consider the following cases.

Case 1. $s = y$. Then $y\in nil(M)$, which is a contradiction.

Case 2. $s = 0$. Then $r_1t = t$, that is, $(r_1 - 1)t = 0$, which implies  $r_1 - 1 \in ann(t)$ and therefore, $M = Rt \oplus ann(t)M$.

Case 3. $s = m$. Then $r_1t - t \in nil(M)$. Let $r_2 = r_1^{2} - r_1$. Therefore $r_2\in [r_1t - t : M]$ and by Lemma 3.7, we have $r_2^{n}(r_1t - t) = 0$, for some $n\in \mathbb{N}$. Thus $r_2^{n + 1}t = 0$. For some suitable choice $s_2$ in terms of $r_2$ we see that $(s_2^{2} - s_2)(1 + 4r_2)t + r_2t  = 0$. Therefore for $v = r_1 + s_2(1 - 2r_1)$, we have $vt = v^{2}t$ and $v\in [t : M]$. It follows that there is some $w\in R$ such that $w\in ann(r_1t)$. Thus by a similar argument as in case 2 above we conclude that $M = Rr_1t \oplus ann(r_1t)M$.\qed\\

\noindent{\bf Remark 3.9.} Theorem 3.8 is also true for the vertices of annihilating graphs $ann_s(\Gamma(M))$ and $ann_t(\Gamma(M))$. That is, if we have the information regarding the elements of objects $\widehat{A_s(M)}$ and $\widehat{A_t(M)}$, we can establish the structure of a module $M$.\\

Recall that a complete bipartite graph is one whose each vertex of one partite set is joined to every vertex of the other partite set. We denote the complete bipartite graph with partite sets of size $m, n\in\mathbb{N}$ by $K_{m,n}$. More generally a complete $r$-partite graph is one whose vertex set can be partitioned into $r$ subsets so that no edge has both ends in any one subset and each vertex of a partite set is joined to every vertex of the another partite sets. A complete bipartite graph of the form $K_{1,n}$ is called a star graph.\\

In the following result, we discuss the neighbourhood relation for the elements $x, y\in M$ with $[x : M][y : M]M \neq 0$ in $ann_f(\Gamma(M))$. We show that the neighbourhood similar vertices determines the nature of the annihilating graphs arising from $M$.\\

\noindent{\bf Theorem 3.10.} Let $M$ be an $R$-module and let $ann_f(\Gamma(M))$ be the full annihilating graph of $M$ which contains no cycle of length 3. If for every pair $z, y\in \widehat{A_f(M)}$ with $[z : M][y : M]M \neq 0$ there exists $x\in \widehat{A_f(M)}\setminus\{y, z\}$ such that $y \sim_{nbd} z$ and $[x : M][y : M]M = 0 = [x : M][z : M]M$, then $ann_f(\Gamma(M))$ is a star graph.

\noindent{\bf Proof.} By Lemma 3.7 (i), we have either, $8 \leq |M| \leq 16$ or $|M| \geq 17$ and $nil(M) = \{0, x\}$ with $2x = 0$, $0 \neq x\in M$. If $8 \leq |M| \leq 16$, then again by Lemma 3.7 (i), $ann_f(\Gamma(M))$ is star graph with at most 5 edges. 

Claim 1: $ann_f(\Gamma(M))$ is an infinite graph, for $|M| \geq 17$, with $nil(M) = \{0, x\}$. 

Let $y\in \widehat{A_f(M)}$ such that $[x : M][y : M]M = 0$. Then for all $r\in [y : M]$, we have $[x : M][ry : M]M = 0$, where $r = \sum\limits_{i = 1}^{n}r_i$, $r_i\in [y : M]$ with $1\leq i \leq n$. If $z\in \widehat{A_f(M)}$ is any other element such that $[x : M][z : M]M = 0 = [ry : M][z : M]M$. Then $ann_f(\Gamma(M))$ contains a cycle on $x$, $z$ and $ry$, a contradiction. Similarly, it can be shown that for all $p$, with $1\leq p\leq m$, $[x : M][r^py : M]M = 0$ and consequently by Lemma 3.7 (ii), $r^py$ is a pendant vertex in $ann_f(\Gamma(M))$. Note, that all $r^py$ are distinct. If not, then $r^py = r^qy$, with $1\leq p < q \leq m$. Therefore $y(r^p - r^q) = yr^{p}(1 - r^{q - p}) = 0$. This implies $(1 - r^{q - p})\in ann(r^py)$.  By Lemma 3.7 (i), $x = \sum\limits_{finite}sm = s^{*}m^{*}$ for some $s^{*}\in [x : M]$ and  $m^{*}\in M$. Thus, $(1 - r^{q - p})m^{*}\in ann(r^py)M = \{0, x\}$. So, either $m^{*} - r^{q - p}m^{*} = 0$ or $m^{*} - r^{q - p}m^{*} = x$. If $m^{*} = r^{q - p}m^{*}$, then
\begin{equation*}
x = s^{*}m^{*} = s^{*}r^{q - p}m^{*} \in s^{*}r^{q - p - 1})m^{*}Rr \subseteq [x : M][r^{q - p - 1}y : M]M = 0, 
\end{equation*} 
which is a contradiction. If $m^{*} - r^{q - p}m^{*} = x$, then 
\begin{equation*}
x - r^{q - p}s^{*}m^{*} = s^{*}m^{*} - r^{q - p}s^{*}m^{*} = s^{*}x \subseteq [x : M]^{2}M = 0,
\end{equation*}
again a contradiction. Thus all $r^py$ are distinct and consequently, $ann_f(\Gamma(M))$ is an infinite graph. 

Claim 2: $ann_f(\Gamma(M))$ is a star graph with center $x$. 

Let $a, b\in \widehat{A_f(M)}$ such that $[a : M][b : M]M \neq 0$ and $[x : M][a : M]M = 0 = [x : M][b : M]M$. Suppose that $ann_f(\Gamma(M))$ is not a star graph. There is some $c \in \widehat{A_f(M)}\setminus\{x, a\}$ such that $[a : M][c : M]M = 0$, that is $a$ is not a pendant vertex in $ann_f(\Gamma(M))$. Consider an element $tc$, where $t = \sum\limits_{\alpha\in \Lambda} t_{\alpha}$, $t_{\alpha}\in [a : M]$ for each $\alpha\in \Lambda$, $\Lambda$ is any finite index set. We show that $tc\not\in \{0, a, x, c, b\}$. If $ tc = 0$, then for $a\in \widehat{A_f(M)}$ we have $x, c\in \widehat{A_f(M)}$ with $x \sim_{nbd} c$. Therefore by Lemma 3.7 (ii), $a$ is a pendant vertex in $ann_f(\Gamma(M))$, which contradicts our supposition. If $tc = x$, then $[x : M][tc : M]M = 0$. Therefore $c\in ann(at)M = \{0, x\}$, a contradiction. If $tc = c$, then $[x : M][tc : M]M = 0$, a contradiction, since $ann_f(\Gamma(M))$ contains no cycle of length 3. Finally, if $tc = b$, then $[c : M][tc : M]M = 0$, which implies, $c\in nil(M) = \{0, x\}$, again a contradiction. Thus $tc\in \widehat{A_f(M)}\setminus\{x, a, c, b\}$. So, there is some $a_1\in \widehat{A_f(M)}$ such that $[a_1 : M][tc : M]M = 0$. It is easy to verify that $a_1\not\in \{0, tc, a, x, c, b\}$. Moreover, it is clear that $[x : M][a_1 : M]M \neq 0$, otherwise $ann_f(\Gamma(M))$ contains cycle of length 3. Let $a_1 = \lambda m_1$, where $\lambda = \sum\limits_{\beta\in \Lambda} \lambda_{\beta}$, $\lambda_{\beta}\in [a_1 : M]$ for each $\beta\in \Lambda$ and $m_1\in M$. Then $[\lambda x : M][a : M]M = 0$, which implies $\lambda x\in ann(a)M$. Thus $\lambda x = x$, note that $\lambda x \neq 0$, otherwise $[x : M][a_1 : M]M = 0$. On the other hand $[a_1 : M][tc : M]M = 0$. Therefore, $[\lambda c : M][a : M]M = 0 = [ta_1 : M][c : M]M$, which implies $\lambda c\in ann(a)M$. Thus $\lambda c = 0$ or $\lambda c = x$. So, $[x : M][c : M]M = 0 = [x : M][b : M]M$, a contradiction. Thus we have a contradiction in every possible case and hence it follows that $ann_f(\Gamma(M))$ is a star graph with center $x$.\qed\\

We conclude this section with the following remark.\\

\noindent{\bf Remark 3.11.} Theorem 3.10 is not true for general simple connected graphs. Consider the graph $\Gamma$ with vertex set $V = \{a_1, a_2, a_3, a_4\}$ and edge set $E = \{a_1-a_2, ~a_2-a_3, ~a_3-a_4, ~a_4-a_1\}$. Clearly, $a_1 \sim_{nbd} a_3$ with $N(a_1) = N(a_3) = \{a_2, a_4\}$, but $\Gamma$ is not a star graph in fact $\Gamma$ is a cycle graph on four vertices. \\

\section{\bf Graph isomorphism and graphs arising from tensor product}

In this section, we discuss isomorphism of annihilating graphs. We exhibit certain module theoretic properties of modules $M$ and $N$ which they share if their full annihilating graphs $ann_f(\Gamma(M))$ and $ann_f(\Gamma(N))$ are isomorphic. Moreover, we consider the annihilating graph arising from tensor product $M\otimes_R T^{-1}R$, where $T = R \backslash C(M)$, where $C(M) = \{ r\in R$ : $rm = 0~ for~ some~ 0\neq m \in M$\}. We investigate the case when $M$ is a multiplication module and show that $ann_f(\Gamma(M)) \cong ann_f(M\otimes_R T^{-1}R)$ for every module $M$.\\

\noindent{\bf Definition 4.1.} Let $x\in M$ be a vertex in $ann_f(\Gamma(M))$. We say that $x$ is primitive vertex if the submodule  generated by $x$ is cyclic over $R$.\\

It can be easily checked that an element $a\in Z^{*}(R)$, where $R$ is a von Neumann regular ring is a primitive vertex in the zero-divisor graph \cite{AdLs} if and only if the ring $Ra$ is a field. We define the order of vertex $x\in\widehat{A_f(M)}$ by $\Theta(x) = |Rx|$. Clearly, if $x\sim_{nbd}y$ for $y\in\widehat{A_f(M)}$, then $|Rx| = |Ry|$, which implies  $\Theta(x) = \Theta(y)$. Thus we can talk about the order of equivalence classes.\\

The following result is one of the entanglement for modules $M$ and $N$ if their full annihilating graphs are isomorphic.\\

\noindent{\bf Theorem 4.2.} Let $M$ and $N$ be two $R$-modules such that $ann_f(\Gamma(M)) \cong ann_f(\Gamma(N))$. If $Soc(M)$ is a sum of finite simple cyclic submodules, then $Soc(M) \cong Soc(N)$.

\noindent{\bf Proof.} Let $\eta : ann_f(\Gamma(M)) \longrightarrow ann_f(\Gamma(N))$ be an isomorphism of graphs. It is clear that if $x\sim_{nbd}y$, then $\eta(x)\sim_{nbd}\eta(y)$ and $\Theta(x) = \Theta(\eta(x))$. Also, the primitive elements of $ann_f(\Gamma(M))$ are in bijection with the primitive elements of $ann_f(\Gamma(N))$. Thus if $x$ is primitive, so is $\Theta(x)$. Suppose that $Soc(M)$ contains the sum of simple cyclic submodules $\delta$ times. Then there are $\delta$ many equivalence classes of primitive vertices in $ann_f(\Gamma(M))$. This implies that $ann_f(\Gamma(N))$ must have $\delta$ equivalence classes of primitive vertices. It follows that the same number of copies of each simple cyclic submodule must occur in both $Soc(M)$ and $Soc(N)$. Thus, $Soc(M)$ and $Soc(N)$ are isomorphic.\qed\\ 

\noindent{\bf Corollary 4.3.} Let $M =\prod\limits_{i\in I}M_i$ and $N = \prod\limits_{i\in I}N_i$, where $M_i$, $N_i$ are finite simple cyclic modules for all $i\in I$ and $I$ is an index set. If $ann_f(\Gamma(M)) \cong ann_f(\Gamma(N))$, then $M \cong N$.

\noindent{\bf Proof.} $Soc(M) = \sum\limits_{i\in I}M_i$ and $Soc(N) = \sum\limits_{i\in I}N_i$. Thus the result follows.\qed\\

\noindent{\bf Corollary 4.4.} Let $M$ and $N$ be two $R$-modules such that $ann_f(\Gamma(M)) \cong ann_f(\Gamma(N))$. If $M$ has an essential socle, then so does $N$.\\

Let $\mathbb{Z}_m \otimes \mathbb{Z}_n$ be tensor product of two finite abelian groups. It is easy to verify that if g.c.d of $m, n\in \mathbb{Z}$ is 1, then $\mathbb{Z}_m \otimes \mathbb{Z}_n = 0$ and in general $\mathbb{Z}_m \otimes \mathbb{Z}_n \cong \mathbb{Z}_d$, where $d$ is g.c.d of $m$ and $n$. It follows that if g.c.d of $m$ and $n$ is 1, then $\reallywidehat{A_f(\mathbb{Z}_m \otimes \mathbb{Z}_n}) = \phi$. However, if g.c.d of $m$ and $n$ is $d$, $d > 1$ and $\mathbb{Z}_d$ is not a simple finite abelian group, then $\reallywidehat{A_f(\mathbb{Z}_m \otimes \mathbb{Z}_n}) \neq \phi$, in fact the graphs $ann_f(\Gamma(\mathbb{Z}_m \otimes \mathbb{Z}_n))$ and $ann_f(\Gamma(\mathbb{Z}_d))$ are isomorphic. Furthermore, if $\mathbb{Z}_p$, $\mathbb{Z}_q$ and $\mathbb{Z}_r$ are any three finite simple abelian groups, where $p, q, r\in \mathbb{Z}$ are primes, then we have the following equality between the combinatorial objects,
\begin{center}
$\reallywidehat{A_f(\mathbb{Z}_p \oplus \mathbb{Z}_q \otimes \mathbb{Z}_p \oplus \mathbb{Z}_r)} = \reallywidehat{A_f(\mathbb{Z}_p\oplus \mathbb{Z}_r)}$.
\end{center} 
It follows that the full-annihilating graph arising from the tensor product $\mathbb{Z}_p \oplus \mathbb{Z}_q \otimes \mathbb{Z}_p \oplus \mathbb{Z}_r$ with vertex set $\reallywidehat{A_f(\mathbb{Z}_p \oplus \mathbb{Z}_q \otimes \mathbb{Z}_p \oplus \mathbb{Z}_r)}$ is same as the full-annihilating graph arising from the direct sum $\mathbb{Z}_p\oplus \mathbb{Z}_r$.\\

Now, we study the annihilating graph structures arising from the tensor product $M \otimes_{R} T^{-1}R$. The following result proved in \cite{ALS} is perhaps the first result which establishes a connection between the graph structure of $R$ and its localization $S^{-1}R$ (total quotient ring) at $S$, where $S = R\setminus Z(R)$.\\

\noindent{\bf Theorem 4.5.} Let $R$ be a commutative ring with unity and let $S^{-1}R$ be the localization of $R$ at $S$. Then the graphs $\Gamma(R)$ and $\Gamma(S^{-1}R)$ are isomorphic.\\

The equivalence class corresponding to equivalence relation $\sim_{M}$ for $m\otimes{ 1 \over {s}}\in M\otimes_{R} T^{-1}R$ is denoted by,
\begin{center}
$[m\otimes{ 1 \over {s}}]_{M\otimes_{R} T^{-1}R} = \{m^{*}\otimes{ 1 \over {s}}\in M\otimes_{R} T^{-1}R : m^{*}\sim_{M}m\}$.
\end{center}
For any module $M$, it can be easily seen that the localization $T^{-1}M$ at $T$ is the special case of tensor product with $T^{-1}M \cong M \otimes_{R} T^{-1}R$. We first investigate the cardinalities of equivalence classes corresponding to the equivalence relation $\sim_M$ on sets $\widehat{A_f(M\otimes_{R} T^{-1}R)}$ and $\widehat{A_f(M)}$.\\ 

\noindent{\bf Lemma 4.6.} Let $M$ be an $R$-module. Then the equivalence classes corresponding to equivalence relation $\sim_M$ on sets $\widehat{A_f(M\otimes_{R} T^{-1}R)}$ and $\widehat{A_f(M)}$ have the same cardinality.

\noindent{\bf Proof.} It suffices to prove the result for sets $\widehat{A_f(M\otimes_{R} T^{-1}R)}$ and $\widehat{A_f(M)}$. First note that for all $m\in\widehat{A_f(M)}$ and for all modules $M, N$ we have,
\begin{center}
$ann_{T^{-1}R}(m\otimes {1\over{s}}) = ann(m)\otimes T^{-1}R$ \hskip 0.4cm and
\end{center}
\begin{center}
\hskip 0.4cm $[N : M] \otimes T^{-1}R ~M\otimes T^{-1}R  = [N \otimes T^{-1}R : M\otimes T^{-1}R]M\otimes T^{-1}R$. 
\end{center}
Moreover,
\begin{center}
$\widehat{A_f(M\otimes T^{-1}R)} = \widehat{A_f(M)}\otimes T^{-1}R$  \hskip 0.5cm and
\end{center}
\begin{center}
 $[m]_M\otimes T^{-1}R = [{m\over 1}]_{M\otimes T^{-1}R}\otimes T^{-1}R$.
\end{center}
This implies,
\begin{center}
 $\widehat{A_f(M)} = \bigcup\limits_{\lambda\in\Lambda}[m_\lambda]_M$ and
\end{center}
\begin{center}
\hskip 0.4cm$\widehat{A_f(M\otimes T^{-1}R)} =  \bigcup\limits_{\lambda\in\Lambda}[m_\lambda\otimes {1\over1}] = \bigcup\limits_{\lambda\in\Lambda}[{m_\lambda \over 1}]_{M\otimes T^{-1}R}$.
\end{center}
Here, $\Lambda$ is an index set with $[m_\lambda] = [m_\mu] = \emptyset$, for $\lambda, \mu \in \Lambda$. We show that for all $m\in\widehat{A_f(M)}$,
$\mid[{m\over 1}]_{M\otimes T^{-1}R}\mid ~= ~\mid[m]_M\mid$.

We consider the following cases.

\noindent Case 1. $[x : M] \neq ann(M)$. Then by [Lemma 3.2, \cite{SR}], $M$ is a multiplication module.
 
\noindent Subcase 1.1. $\mid [m]_M\mid$ is finite.

Let ${m^{*}\over s} \in [m\otimes {1\over s}]_{M\otimes T^{-1}R}$ with $m^{*}\in[m]_M$, $s\in T$. Then $ann(m^{*})M = ann(m)M$ and $\{s^{n}m^{*} : n\geq 1\} \subseteq [m]_M$. Therefore,

\begin{center}
$m^{*}\otimes {1\over s} = {m^{*}\over s} = {m^{*}s^{k}\over s^{k+1}} = m^{*} \in [m]_M$.
\end{center} Since $\mid [m]_M \mid$ is finite, so there exists some $k\in \mathbb{N}$ such that $s^{k}m^{*} = s^{k+1}m^{*}$. Thus $\widehat{A_f(M\otimes T^{-1}R)} \subseteq \widehat{A_f(M)}$ which implies $[{m\over 1}]_{M\otimes T^{-1}R} \subseteq [m]_M$ and we have,
\begin{center}
$\mid[{m\over 1}]_{M\otimes T^{-1}R}\mid ~= ~\mid[m]_M\mid$.
\end{center} 
Subcase 1.2. $\mid [m]_M\mid$ is infinite.

\noindent We define an equivalence relation $\sim$ on $T$ by $s \sim x$ if and only if $sm = xm$. Define a map 
\begin{center}
$\varphi : M \times T/\sim ~\longrightarrow [{m\over 1}]_{M\otimes T^{-1}R}$
\end{center}
\begin{center} 
\hskip 0.2 cm $(a,[s]) \longrightarrow a\otimes{1\over s}$                                           
\end{center}
We first check that the map $\varphi$ is well defined. Let $m_1 = m_2$ and $[s] = [x]$. Then,
\begin{center}
$(s-x)M \subseteq ann(m)M = ann(m_1)M = ann(m_2)M$.
\end{center}
Since $M$ is a multiplication module therefore we have,
\begin{center}
$(s-x)M = ann(m)M = ann(m_1)M = ann(m_2)M$.
\end{center}
Thus $sm_1 = xm_1$ and $sm_2 = xm_2$ which implies,
\begin{center}
${m_1\over s} = {m_2\over x}$.
\end{center}
It is clear that the map $\varphi$ is surjective. Thus,
\begin{center}
$\mid[{m\over 1}]_{M\otimes T^{-1}R}\mid ~\geq~ \mid[m]_M\mid \mid T/\sim \mid$.
\end{center}
Further, the map
\begin{center}
$T/\sim ~\longrightarrow [m]_M$
\end{center}
\begin{center}
\hskip 0.5 cm $[s] ~\longrightarrow sx$ 
\end{center}
is well defined and injective. Therefore, 
\begin{center}
$\mid[{m\over 1}]_{M\otimes T^{-1}R}\mid ~\leq \mid[m]_M\mid^{2} ~=  ~\mid[m]_M\mid$.
\end{center}
Since $\mid[m]_M\mid$ is infinite we conclude that $\mid[{m\over 1}]_{M\otimes T^{-1}R}\mid ~= ~\mid[m]_M\mid$.

\noindent Case 2. $[x : M] = ann(M)$, for some $0 \neq x \in M$. Then for all $0 \neq y\in M$, we have $[x : M][y : M]M = 0$. Thus $\widehat{A_f(M)} = \widehat{M}$ and $\widehat{A_f(M\otimes T^{-1}R)} = \widehat{M\otimes T^{-1}R}$. Using above subcases 1.1 and 1.2, it can be shown that for each $\alpha\in \Lambda$, the cardinality of equivalence classes $[m_{\alpha}]_M$ and $[{m_{\alpha}\over 1}]_{M\otimes T^{-1}R}$ are same.  
\qed\\

\noindent{\bf Theorem 4.7.} Let $M$ be an $R$-module. Then $ann_f(\Gamma(M\otimes T^{-1}R)) \cong ann_f(\Gamma(M))$.

\noindent{\bf Proof.} For each $\alpha\in \Lambda$, $\Lambda$ an index set, we see by Lemma 4.6 that there is a bijection
\begin{center} 
$\varphi_{\alpha}:[m_{\alpha}]_M \longrightarrow [{m_{\alpha}\over 1}]_{M\otimes T^{-1}R}$
\end{center}
Define a map,
\begin{center}
\hskip 0.2cm $\psi :  \widehat{A_f(M)} \longrightarrow \widehat{A_f(M\otimes_{R} T^{-1}R)}$
\end{center}
\begin{center}
$ m \longrightarrow [m_{\alpha}]_M$ 
\end{center}
Clearly $\psi$ is a bijective map. To show the required isomorphism of graphs we just need to show that $m_{1}^{*}$ and $m_{2}^{*}$ are adjacent in $ann_f(\Gamma(M))$ if and only if $\psi(m_{1}^{*})$ and $\psi(m_{2}^{*})$ are adjacent in $ann_f(\Gamma(M\otimes T^{-1}R))$. That is,
\begin{center}
\hskip 0.9cm $[m_{1}^{*} : M][m_{2}^{*}: M]M = 0$  \hskip 0.2cm if and only if 
\end{center}
\begin{center}
$[\psi(m_{1}^{*}) : M][\psi(m_{2}^{*}): M]M = 0$.
\end{center}
It suffices to show that for 
$m_{1}^{*}\in[m]_M$, $m_{1}^{*}\in[r]_M$, $x_1 \in [{m\over 1}]_{M\otimes_{R} T^{-1}R}$ and  $x_2 \in [{r\over 1}]_{M\otimes_{R} T^{-1}R}$
\begin{center}
$[m_{1}^{*} : M][m_{2}^{*}: M]M = 0$ \hskip 0.4cm if and only if 
\end{center} 
\begin{center}
$[{x_1\over 1} : M\otimes T^{-1}R][{x_2\over 1} : M\otimes T^{-1}R]M\otimes T^{-1}R = 0$.
\end{center}

Therefore,

\hskip 0.5 cm $[m_{1}^{*} : M][m_{2}^{*}: M]M = 0$ 

$\iff m_{1}^{*}\in ann(m_{2}^{*})M = ann(r)M$

$\iff {m_{1}^{*}\over 1}\in ann({m_{2}^{*}\over 1})M\otimes T^{-1}R = ann({r\over 1})M\otimes T^{-1}R = ann({x_{2}\over 1})M\otimes T^{-1}R$

$\iff [{m_{1}^{*}\over 1} : M\otimes T^{-1}R][{x_{2}\over 1} : M\otimes T^{-1}R]M\otimes T^{-1}R = 0$

$\iff {x_{2}^{*}\over 1}\in ann({m_{1}^{*}\over 1})M\otimes T^{-1}R = ann({m\over 1})M\otimes T^{-1}R = ann({x_{1}\over 1})M\otimes T^{-1}R$

$\iff [{x_{2}\over 1} : M\otimes T^{-1}R][{x_{1}\over 1} : M\otimes T^{-1}R]M\otimes T^{-1}R = 0$.

Thus it follows that,
\begin{center} $ann_f(\Gamma(M\otimes T^{-1}R)) \cong ann_f(\Gamma(M))$.\end{center}\qed\\

 The following result is an immediate consequence of preceding theorem.\\
 
\noindent{\bf Corollary 4.8.} Let $M$ be an $R$-module. Then the following  hold.

(i) $ann_s(\Gamma(M\otimes T^{-1}R) \cong ann_s(\Gamma(M))$

(ii) $ann_t(\Gamma(M\otimes T^{-1}R) \cong ann_t(\Gamma(M))$\\

\noindent{\bf Remark 4.9.} Consider a $\mathbb{Z}$-module $\mathbb{Z}_2 \oplus \mathbb{Z}_2$. It is easy to verify that
\begin{center} $ann_f(\Gamma(\mathbb{Z}_2 \oplus \mathbb{Z}_2)) \cong ann_s(\Gamma(\mathbb{Z}_2 \oplus \mathbb{Z}_2)) \cong K_3$,
\end{center}
By Theorem 4.7, we see that,
\begin{center} $ann_f(\Gamma(\mathbb{Z}_2 \oplus \mathbb{Z}_2 \otimes T^{-1} \mathbb{Z})) \cong ann_f(\Gamma(\mathbb{Z}_2 \oplus \mathbb{Z}_2)) \cong K_3$,
\end{center}
and
\begin{center} $ann_s(\Gamma(\mathbb{Z}_2 \oplus \mathbb{Z}_2 \otimes T^{-1} \mathbb{Z})) \cong ann_s(\Gamma(\mathbb{Z}_2 \oplus \mathbb{Z}_2)) \cong K_3$,
\end{center}
where as,
\begin{center} $ann_t(\Gamma(\mathbb{Z}_2 \oplus \mathbb{Z}_2 \otimes T^{-1} \mathbb{Z})) \cong ann_t(\Gamma(\mathbb{Z}_2 \oplus \mathbb{Z}_2)) \cong \emptyset$.
\end{center}
It follows that annihilating graphs $ann_t(\Gamma(\mathbb{Z}_2 \oplus \mathbb{Z}_2 \otimes T^{-1} \mathbb{Z}))$ and  $ann_f(\Gamma(\mathbb{Z}_2 \oplus \mathbb{Z}_2 \otimes T^{-1} \mathbb{Z}))$ are not isomorphic. That is, all annihilating graphs for a same tensor product cannot be similar. However, if $M$ is multiplication $R$-module, then by [Theorem 3.9, \cite{SR}] all three annihilating graphs arising from $M$ are same. Therefore by Theorem 4.7, all three annihilating graphs arising from the tensor product $M \otimes T^{-1}R$ are same. In fact, the three annihilating graphs arising from $M \otimes T^{-1}R$ coincides with the three annihilating graphs arising from $M$.\\

We conclude this paper with some discussion on factor modules $M = \prod\limits_{i\in I}M_i/\sum\limits_{i\in I}M_i$ and $N = \prod\limits_{i\in I}N_i/\sum\limits_{i\in I}N_i$, where $I$ is an index set, $M_i$,  $N_i$ are finite simple modules not equal to $\mathbb{Z}_2(\mathbb{Z})$ for any $i\in I$. Consider the equivalence relation on subsets of $I$ given by $J$ is equivalent to $K$ if symmetric difference of $J$ and $K$ is finite. Pick one element $J$ from each equivalence class and for each such set $J$, let $X_J = \{x + \sum\limits_{i\in I}M_i : x(j) = 0 \iff j\in J\} \subset M$ and $Y_J = \{y + \sum\limits_{i\in I}N_i : y(j) = 0 \iff j\in J\} \subset N$. Since a finite simple module $\mathbb{Z}_2(\mathbb{Z})$ is a component of neither $M$ nor $N$, therefore sets $X_J$ and $Y_J$ have same cardinality. Let $\varphi_J$ be some bijection between sets $X_J$ and $Y_J$. It is clear that the graphs $ann_f(\Gamma(M))$, $ann_f(\Gamma(N))$ are complete with $\widehat{A_f(M)} = \widehat{M}$, $\widehat{A_f(N)} = \widehat{N}$ and for each $a\in M$, we have the associated zero set $J$. Thus the map $a$ to $\varphi_J(a)$ defines a bijection between the vertices of $ann_f(\Gamma(M))$, $ann_f(\Gamma(N))$ and we conclude that $ann_f(\Gamma(M))\cong ann_f(\Gamma(N))$.

\end{document}